\documentclass[11pt]{article}
\usepackage[a4paper, margin=0.75in, bmargin=0.85in]{geometry}
\usepackage[mathscr]{euscript}
\usepackage{amsmath,amsthm,mathrsfs}
\usepackage{amsfonts,amssymb}
\usepackage{graphicx}
\usepackage{breqn}
\usepackage{mathrsfs}
\usepackage{amsmath}
\usepackage{float}
\usepackage{caption}
\usepackage{chngcntr}
\usepackage{tikz-qtree}
\usetikzlibrary{arrows.meta}
\usepackage{subcaption}
\usepackage{epsfig,amsmath,amsthm,latexsym,amssymb,amsgen,graphicx}
\DeclareMathAlphabet{\mathpzc}{OT1}{pzc}{m}{it}
\usepackage{verbatim}
\usepackage{cite}
\newtheorem{thm}{Theorem}[section]

\newtheorem{cor}[thm]{Corollary}

\newtheorem{prop}[thm]{Proposition}

\newtheorem{defn}[thm]{Definition}
\newtheorem{rem}[thm]{Remark}

\newtheorem{exam}[thm]{Example}
\begin{document}
\title {\bf{Domination Index in Graphs}}
\author { \small Kavya. R. Nair \footnote{Corresponding author. Kavya. R. Nair, Department of Mathematics, National Institute of Technology, Calicut, Kerala, India. Email: kavyarnair@gmail.com} and M. S. Sunitha \footnote{sunitha@nitc.ac.in} \\ \small Department of Mathematics, National Institute of Technology, Calicut, Kerala, India-673601}
\date{}
\maketitle
\hrule
\begin{abstract}
The concepts of domination and topological index hold great significance within the realm of graph theory. Therefore, it is pertinent to merge these concepts to derive the domination index of a graph. A novel concept of the domination index is introduced, which utilizes the domination degree of a vertex. The domination degree of a vertex $a$ is defined as the minimum cardinality of a minimal dominating set that includes $a.$ The idea of domination degree and domination index is conducted of graphs like complete graphs, complete bipartite, $r-$ partite graphs, cycles, wheels, paths, book graphs, windmill graphs, Kragujevac trees. The study is extended to operation in graphs. Inequalities involving domination degree and already established graph parameters are discussed. An application of domination degree is discussed in facility allocation in a city. Algorithm to find a MDS containing a particular vertex is also discussed in the study. 
\end{abstract}
\hrule
\section{Introduction}
Graphs are used to model pairwise relationships between objects. Graphs consist of vertices (also known as nodes) and edges that connect pairs of vertices.  Graph theory is an inevitable area of research because of its ability to depict and examine relationships, its universality, fundamental concepts, and widespread application in diverse domains. Graph theory has applications in the fields of computer science, operations research, social sciences, biology, and many more. Topological indices are numerical values that quantify a graph's structural traits or attributes. They are used in a variety of disciplines, including chemistry, network analysis, and algorithm creation. These indices offer insights, prediction power, and quantitative representations of graph topologies by capturing essential characteristics of graph connectivity or degrees. Among all the topological indices the most extensively studied and used index is the Wiener index \cite{W} which was developed by Harold Wiener in 1947 to compare the boiling point of some alkane isomers. Due to its widespread application numerous other topological indices were introduced which includes Randic Index (1975) \cite{r}, Zagreb Indices (1972) \cite{z1,z2}, Estrada Index (2000) \cite{e}, Gutman Index (1994) \cite{g}, Harary index (1993) \cite{h1,h2}, Balaban Index (1982) \cite{b1,b2} etc.  Domination in graphs is a fundamental concept that is vital due to its significance in wireless sensor networks, social networks, decision-making, and connection with other graph concepts. Ore and Berge's \cite{ore1962, berge1962} work on graph domination and related topics significantly contributed to the development and understanding of domination theory. The notion of domination topological index in graph was introduced by Hanan et. al \cite{T5} in 2021.
\\Motivated from this idea this article introduces a novel idea of a domination index using minimal dominating sets.
Consider a graph and a vertex $v.$ The domination index is defined using a dominating set containing $v$ having the least cardinality. In a graph the least cardinality of a dominating set is obtained by considering the minimum dominating set. But it is not always possible to find a minimum dominating set containing $v.$ However it is always possible to find a dominating set that is minimal containing $v.$ Hence, in this article the notion of domination degree of a vertex $v$  is defined using the number of vertices in a minimal dominating set containing $v.$ The article is structured as follows:\\ The initial section provides a concise overview of the fundamental of graph theory. Section 3 begins with an important result which states that it is always possible to find a MDS containing a particular vertex in a graph. the section introduces and illustrates the novel concept of domination degree. Inequalities relating domination degree and existing graph parameters are obtained. The idea of domination regular graphs is discussed. The concept of domination degree is studied in various graphs like complete graphs, complete bipartite and $r-$ partite graphs, cycle, paths, book graphs, windmill graphs, wheels, Kragujevac trees. The study is also extended to operations on graphs like union, join, composition and corona products. Section 4 introduces the concept of domination index in graphs. Using the results in Section 2, domination index is studied in various graphs and graph operations. An application of domination degree in facility allocation is provided in Section 5. Section 5 also deals with an algorithm to find a MDS containing a particular vertex.  
\section{Preliminaries}
The following fundamental definitions in graph theory is referred from \cite{W, T5, bal, char, k}
A graph $\mathcal{G}=(\mathcal{V},\mathcal{E})$, consist a set $\mathcal{V}$ of vertices which is nonempty and a set $\mathcal{E}$ of edges which is a 2- element subset of $\mathcal{V}.$ An edge $e=ab$ means that $a$ and $b$ are adjacent or neighbors, and $a,b$ are end vertices of $e.$ Let $a$ be a vertex in $\mathcal{G}.$ The open neighborhood of $a$ is $\mathcal{N}(a)=\{b\,:\, \text{ab is an edge}\}.$ The closed neighborhood of $a$ is $\mathcal{N}[a]=\{a\}\cup \mathcal{N}(a).$ The order $n$ and size $m$ are the number of vertices and edges respectively. Let $\mathcal{G}$ and $\mathcal{H}$ be two simple graphs. Then a graph isomorphism between $\mathcal{G}$ and $\mathcal{H}$ is a bijection $\theta:\mathcal{V}(\mathcal{G})\rightarrow \mathcal{V}(\mathcal{H})$ such that $a$ and $b$ are adjacent in $\mathcal{G}$ iff $\theta(a)$ and $\theta(b)$ are adjacent in $\mathcal{H}.$ A graph $\mathcal{H}$ is a subgraph of $\mathcal{G}$ if $\mathcal{V}(\mathcal{H})\subseteq \mathcal{V}(\mathcal{G})$ and $\mathcal{E}(\mathcal{H})\subseteq \mathcal{E}(\mathcal{H}).$ If $\mathcal{V}(\mathcal{G})=\mathcal{V}(\mathcal{H}),$ then $\mathcal{H}$ is a spanning subgraph. Now, $\mathcal{T}$ is a spanning tree of $\mathcal{G},$ is $\mathcal{T}$ is both a spanning subgraph and a tree.\\
Let $\mathcal{G}_1=(\mathcal{V}_1,\mathcal{E}_1)$ and $\mathcal{G}_2=(\mathcal{V}_2,\mathcal{E}_2)$ be two graphs. Then the union $\mathcal{G}=\mathcal{G}_1\cup \mathcal{G}_2$ is a graph with $\mathcal{V}=\mathcal{V}_1\cup \mathcal{V}_2$ and $\mathcal{E}=\mathcal{E}_1\cup \mathcal{E}_2.$ The join $\mathcal{G}=\mathcal{G}_1+\mathcal{G}_2$ is a graph with $\mathcal{V}=\mathcal{V}_1\cup \mathcal{V}_2$ and $\mathcal{E}=\mathcal{E}_1\cup \mathcal{E}_2\cup \mathcal{E}',$ where $E'$ is the set of edges that join each vertex of $\mathcal{V}_1$ to every vertex of $\mathcal{G}_2.$ Let $\mathcal{G}=\mathcal{G}_1\times\mathcal{G}_2$ be the Cartesian product of $\mathcal{G}_1$ and $\mathcal{G}_2.$ Then $\mathcal{V}(\mathcal{G})=\mathcal{V}_1\times \mathcal{V}_2.$ Two vertices $(a_1,b_1)$ and $(a_2,b_2)$ are adjacent in $\mathcal{G},$ iff either $a_1=a_2$
and $b_1$ is adjacent to $b_2$ in $\mathcal{G}_2$ or $a_1$ is adjacent to $a_2$ in $\mathcal{G}_2$ and $b_1=b_2.$ The vertex set of composition and direct product of graphs are also $\mathcal{V}_1\times \mathcal{V}_2.$ Two vertices $(a_1,b_1)$ and $(a_2,b_2)$ are adjacent in composition $\mathcal{G}= \mathcal{G}_1\circ \mathcal{G}_2,$ iff f $a_1$ is adjacent to $a_2$
in $\mathcal{G}_1$ or $a_1=a_2$ and $b_1$ is adjacent to $b_2$ in $\mathcal{G}_2.$ Similarly, Two vertices $(a_1,b_1)$ and $(a_2,b_2)$ are adjacent in direct product $\mathcal{G}= \mathcal{G}_1\square \mathcal{G}_2,$ iff f $a_1$ is adjacent to $a_2$
in $\mathcal{G}_1$ and $b_1$ is adjacent to $b_2$ in $\mathcal{G}_2.$ The strong product $\mathcal{G}=\mathcal{G}_1\boxtimes \mathcal{G}_2$ is the union of direct product and Cartesian product. The corona $\mathcal{G}=\mathcal{G}_1\odot\mathcal{G}_2$ is obtained by taking $|\mathcal{V}_1|$ copies of $\mathcal{G}_2$ and joining $i^{th}$ vertex if $\mathcal{G}_1$ to every vertex in the $i^{th}$ copy of $\mathcal{G}_2.$  \\
A wheel graph $\mathcal{W}_n$ is the join of $\mathcal{K}_1$ and cycle $\mathcal{C}_n.$ A Windmill graph $\mathcal{W}d(r,s)$ is obtained by taking $s$ copies of $\mathcal{K}_r$ with a common vertex. The Cartesian product of a star graph $\mathcal{S}_{n+1}$ and path $\mathcal{P}_2$ is a book graph $\mathcal{B}_n.$ 
Now, consider $\mathcal{P}_3.$ For $s=2,3,...$ identify the roots of $s$ copies of $\mathcal{P}_3$ to obtain a rooted tree $\mathscr{B}_k$, and the root of $\mathscr{B}$ is the vertex obtained by identifying copies of $\mathcal{P}_3.$
Let $t\geq 2,$ and $\beta_1,\beta_2,...\beta_t\in \{\mathscr{B}_1,\mathscr{B}_2,...\}.$ A Kragujevac tree, $\mathcal{T}$, is a tree that includes a vertex say $a$ with a degree of $t$. The vertex $a$ is adjacent to the roots of $\beta_1,\beta_2,...\beta_t.$ Vertex $a$ is referred to as the central vertex of $\mathcal{T}$. The subgraphs $\beta_1,\beta_2,...\beta_t$ are called the branches of $\mathcal{T}$.
Let $\mathcal{D}\subseteq \mathcal{V},$ then $\mathcal{D}$ is a dominating set (DS) if for every vertex  $b\in \mathcal{V}\setminus \mathcal{S}$ there exists $a\in \mathcal{S}$ such that $a$ dominates $b, $ i.e, $a$ is adjacent to $b.$ The minimum of number of vertices in a DS in $\mathcal{G}$ is called domination number of $\mathcal{G}$ denoted by $\gamma(\mathcal{G})$ or simply $\gamma.$ The DS with least number of vertices is called a minimum DS. A DS is minimal if it does not properly contain any other DS. The private neighborhood of vertex $a$ with respect to set $\mathcal{S}$, denoted as $P_{N}[a,\mathcal{S}]=\mathcal{N}[a]\setminus \bigcup\limits_{b\in \mathcal{S}\setminus\{a\}}\mathcal{N}[b].$ Here, $\mathcal{S}$ is an irredundant set if every vertex in $\mathcal{S}$ has non- empty private neighborhood. Upper domination number of $\mathcal{G}$, denoted as $\Gamma(\mathcal{G})$ or simply $\Gamma,$ is the maximum cardinality of a minimal DS in $\mathcal{G}.$ Similarly, the irredundance number, $ir(\mathcal{G})$ is the minimum cardinality of maximal irredundant set of $\mathcal{G}.$ And the maximum cardinality of an irredundant set is the upper irredundance number, $IR(\mathcal{G}.)$
\begin{thm}\label{MDSch}\cite{bal}
    A DS $\mathcal{D},$ of graph $\mathcal{G}$ is MDS iff $P_{N}[a,\mathcal{D}]\neq \phi,$ $\forall a\in \mathcal{V}(\mathcal{G}).$
\end{thm}
\begin{thm}\label{MDSDS}\cite{bal}
Let $\mathcal{G}$ be a graph without isolated vertices. If $\mathcal{D}$ is a MDS then $\mathcal{V}\setminus \mathcal{D}$ is a DS. 
\end{thm}
\begin{thm}\label{domin}\cite{bal}
    For a graph $\mathcal{G},$ $\lceil{\frac{n}{1+\Delta(\mathcal{G})}}\rceil \leq \gamma(\mathcal{G})\leq n-\Delta(\mathcal{G}).$
\end{thm}
\begin{thm}\label{irin}\cite{bal}
    For any graph $\mathcal{G},$ $ir(\mathcal{G})\leq \gamma(\mathcal{G})\leq \Gamma(\mathcal{G})\leq IR(\mathcal{G}).$
\end{thm}
\section{Domination degree of a vertex}
\begin{figure}
\begin{center}
    \includegraphics[height=4cm,width=3.5cm]{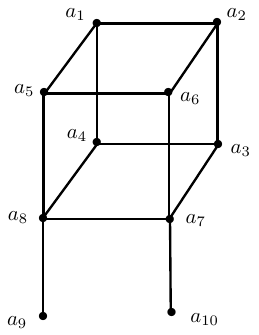}
    \caption{Illustration of domination degree of a vertex.}
    \label{f1}
\end{center}
\end{figure}

The section introduces a new degree of a vertex say $a$ using the minimal dominating sets containing vertex $a$. The inequalities relating various graph parameters and the defined domination degree are obtained. The study is conducted on various graph families and significant graphs. Section begins with the result that ascertain the existence of a MDS containing a particular vertex. The section also defines a domination regular graph (DRG). The study is conducted on graph operations like union, join, composition and corona. The domination degree of vertices of Petersen graph, Herschel graph and Gr\"{o}tzsch graph is also considered.    
\begin{thm}\label{MDS}
Let $\mathcal{G}$ be a connected graph and let $a$ be any vertex in $\mathcal{G}$, then there always exists a MDS containing $a.$     
\end{thm}
\begin{proof}
   Let $\mathcal{G}$ be a connected graph and $a\in \mathcal{V}(\mathcal{G}).$ Suppose $a\in \mathcal{D},$ where $\mathcal{D}$ is a minimum DS, then $\mathcal{D}$ is a MDS containing $a.$
   Suppose that no minimum DS contains $a$ and let $\mathcal{D}'$ be a minimum DS. Then $\mathcal{V}\setminus \mathcal{D}'$ is a DS containing $a.$ Then by a constructive proof, a MDS can be obtained from $\mathcal{V}\setminus \mathcal{D}$ containing $a.$ 
   Now, consider any DS $\mathcal{S}$ containing $a.$ The aim is to obtain a MDS containing $a$ from $\mathcal{S}.$ 
   Suppose $\mathcal{S} $ is not minimal. Then $\exists$ a vertex $b$ such that $\mathcal{S} \setminus \{b\}$ is a dominating set. If $b\neq a$ and $\mathcal{S} \setminus \{b\}$ is MDS then $\mathcal{S} \setminus \{b\}$ is the required MDS. If not repeat the process and in each stage if $b$ is not the deleted vertex then finally a set containing $b$ is obtained which is a MDS.\\
   Now, suppose by deleting $a$, a DS is obtained. This means that $\mathcal{N}[a]$ is dominated by vertices of $\mathcal{S}.$
    \\
    \textbf{Case 1:} $\mathcal{N}[a]$ is dominated by the same vertex say $b$ in $\mathcal{S}$.\\
    Suppose that $N[a]=\{a,a_1,a_2,...,a_n\}$. 
    Also let $N_s[b]\setminus \{a,a_1,...a_n\}=\{b_1,b_2,...,b_m,b_1',b_2',...,b_k'\}$ where $b_1,b_2,..,b_m$ are end vertices. Then $b_',b_2',...b_k'$ are neighbors of $b$ that are not end vertices, i.e, $v$ and $b_i'$ belongs to a common connected component $\forall i=1,2,...,k.$ Also some vertices in $b_1',b_2,...,b_k'$ may be adjacent to each other. Now, delete $b$ from $\mathcal{S}$ and add the following vertices:
    \begin{enumerate}
        \item All the end vertices that are neighbors of $b$, i.e, $\{b_1,b_2,...,b_n\}.$
        \item If two vertices in $\{b_1',b_2',...,b_k'\}$ are neighbors then add any one of them. 
        \item If while considering two vertices $b_i'$ and $b_j'$ in $\{b_1',b_2',...,b_k'\}$ are not neighbors, then add both of them.
    \end{enumerate}
    Let $\mathcal{S}'$ be the set obtained after performing steps 1, 2 and 3.
    As the vertex that dominates $a$ and all its neighbor is deleted $P_N[a,\mathcal{S}']=\phi.$ If after performing steps 1,2 and 3, the minimality condition gets disturbed then delete the vertex with empty private neighborhood. Finally the set obtained is a MDS containing $a. $\\
    \textbf{Case 2:} $N[a]$ is dominated by distinct vertices in $\mathcal{S}$.\\
    Suppose $a$ is dominated by vertex $b$ in $\mathcal{S}$ and the neighbors of $a$ are dominated by $c_1,c_2,..c_{k'}$ in $\mathcal{S}$. 
    Then repeat the same steps $1,2$ and $3$ with the neighbors of $b,c_1,c_2,...c_{k'}.$ and delete the vertex with empty neighborhood. The case when $a$ and its neighbors being dominated by multiple vertices of $\mathcal{S}$ can be proved similarly. \\
    Hence a MDS containing a particular vertex $a$ is obtained from a DS containing $a.$ Therefore, it is always possible to find a MDS containing a particular vertex in a graph $\mathcal{G}$.   
\end{proof}
\begin{defn}\label{dd}
Let $\mathcal{G}$ be a graph and $a\in \mathcal{V}.$ Then the domination degree of vertex $a$ is defined as the minimum number of vertices in the MDS containing $a,$ i.e, 
$$\mathpzc{d}_{\mathpzc{d}}(a)= \min \{|\mathscr{D}| \,: \, \mathscr{D} \text{ is MDS containing $a$ }\}.$$ To be specific the domination degree of a vertex $a$ in $\mathcal{G}$ is denoted as $_{\mathcal{G}} \mathpzc{d}_{\mathpzc{d}}(a).$ Accordingly, the minimum and maximum domination degree of a graph are given by,
$$_{\mathpzc{d}}\delta(G)= min\,\{\mathpzc{d}_{\mathpzc{d}}(a):\, a\in \mathcal{V}(\mathcal{G})\}$$
$$_{\mathpzc{d}}\Delta(G)= max\,\{\mathpzc{d}_{\mathpzc{d}}(a):\, a\in \mathcal{V}(\mathcal{G})\}.$$
\end{defn}
An example to illustrate Definition \ref{dd} is provided in Example \ref{e1}.
\begin{exam}\label{e1}
    Consider the graph in Figure \ref{f1}. Consider vertex $a_1.$ Here, $\mathcal{N}[a_1]=\{a_1, a_2, a_5, a_4\}.$ Hence $a_1$ dominates 4 vertices. To dominate $a_9$ either $a_9$ or $a_8$ is required. Similarly, to dominate $a_{10}$, either $a_{10}$ or $a_7$ is required. Hence a MDS containing $a_1$ with least number of vertices is $\{a_1,a_8,a_{9}\}.$ Therefore, $\mathpzc{d}_{\mathpzc{d}}(a_1)=3.$ Also for vertex $a_5$, 3 any 3 vertices including $a_5$ is insufficient to dominate the whole vertex set. And $\{a_5,a_3,a_9,a_{10}\}$ is a MDS containing $a_5$ with least number of vertices. Hence $\mathpzc{d}_{\mathpzc{d}}(a_5)=4.$ 
\end{exam}
\begin{prop}\label{in1}
Let $\mathcal{G}$ be a graph. Then $\gamma(\mathcal{G})\leq\mathpzc{d}_{\mathpzc{d}}(a)\leq \Gamma(\mathcal{G}),$ $\forall a\in \mathcal{V}(\mathcal{G}).$  
\end{prop}
\begin{proof}
    The proof directly follows from the definition of domination number, domination degree of a vertex [Definition \ref{dd}] and the upper domination number.
\end{proof}
\begin{thm}\label{in2}
 Let $\mathcal{G}$ be a connected graph with $n>1.$ Then, $\mathpzc{d}_{\mathpzc{d}}(a)\leq WI(\mathcal{G}), \forall a \in \mathcal{V}(\mathcal{G}).$     
\end{thm}
\begin{proof}
    Let $\mathcal{G}$ be a connected graph with $n>1.$ Consider a vertex $a$ and let $\mathcal{D}$ be the MDS containing $a$ with the least number of vertices. \\
    \textbf{Case 1}: Suppose $\mathcal{D}=\{a\}.$\\
    Since $n>1$ and $\mathcal{G}$ is connected, there exists at least one edge incident at $a.$ Hence, $\mathpzc{d}_{\mathpzc{d}}(a)\leq WI(\mathcal{G}).$
    \\
    \textbf{Case 2}: Suppose $\mathcal{D}=\{a,b\}.$\\
    In this case, if $a$ and $b$ are adjacent and since $\mathcal{D}$ is minimal, there is at least one vertex each, say $c$, $d$ in the private neighborhood of $a$ and $b$ respectively. Hence there will be at least 3 edges in $\mathcal{G}.$ Similarly, if $a$ and $b$ are non- adjacent, since $\mathcal{G}$ is connected there exists at least one vertex in the path from $a$ to $b$. Hence there are at least two edges in $\mathcal{G}.$ Therefore, $\mathpzc{d}_{\mathpzc{d}}(a)\leq WI(\mathcal{G}).$\\
    \textbf{Case 3}: Suppose $\mathcal{D}=\{a,a_1,a_2,...,a_n\}.$\\
    Each vertex $u\in \{a,a_1,a_2,...,a_n\}$ has either a private neighbor $v$ such that $u\neq v$ or $P_N[u,D]=\{u\}$, which means $u $ is isolated in $\mathcal{G}[D]$. Among $a,a_1,a_2,...,a_n$ let $u_1,u_2,...,u_k$ be the vertices with at least one  distinct vertex in the private neighborhood, i.e, $\exists$ vertices $v_1,v_2,...,v_k$ such that $v_i\in P_N[u_i,D], 1\leq i\leq k .$ And let $x_1,x_2,...,x_{k'}$, $k+k'=n+1$ be the vertices such that $P_N[x_i,D]=\{x_i\}, 1\leq i\leq k'.$ 
     For vertex $u_i$, there exists at least one edge incident at $u_i$ i.e, there exists $u_i v_i$ at $u_i,\forall 1\leq i\leq k.$  Consider two vertices among $x_1,x_2,...,x_{k'}$. A path from $x_i$ to $x_j$, $1\leq i,j\leq k', i\neq j$ contains at least one vertex $x_{i'}$ such that $x_{i'}\notin N_s[u],$ $u\in \{v_1,v_2,...,v_k\} $. Therefore, $\exists$ at least two edges in each path. Suppose there is only one vertex $v$ such that $P_N[v, D]=\{v\}$. Then $v$ is isolated in $\mathcal{G}[D]$. Since $\mathcal{G}$ is connected, consider any path $\mathcal{P}$ from $v$ to $a_i$, then $\exists$ at least two edges in $\mathcal{P}$. 
     Hence, $\mathpzc{d}_{\mathpzc{d}}(a)\leq WI(\mathcal{G}).$
\end{proof}
\begin{rem}\label{dominsub}
It is known result that $\gamma(\mathcal{G})\leq\gamma(\mathcal{G}\setminus e).$ Hence two results follows from this remark:    
\end{rem}
\begin{prop}
    Let $\mathcal{H}$ be a spanning subgraph of $\mathcal{G}.$ Then $_{\mathcal{G}}\mathpzc{d}_{\mathpzc{d}}(a)\leq_{\mathcal{H}}\mathpzc{d}_{\mathpzc{d}}(a)$ for every vertex $a\in \mathcal{V}(\mathcal{G}).$
\end{prop}
\begin{cor}
Let $\mathcal{T}$ be a spanning tree of $\mathcal{G}.$ Then $_{\mathcal{G}}\mathpzc{d}_{\mathpzc{d}}(a)\leq_{\mathcal{T}}\mathpzc{d}_{\mathpzc{d}}(a)$ for every vertex $a\in \mathcal{V}(\mathcal{G}).$    
\end{cor}
Now, the domination index is studied in some significant graphs and graph families. 
\begin{prop}\label{cmp}
For any vertex $a\in \mathcal{K}_n$ of order n, $\mathpzc{d}_{\mathpzc{d}}(a)=1.$
\end{prop}
\begin{proof}
    For $\mathcal{K}_n,$ each vertex dominates every other vertex. Hence $\{a\},$ $\forall\, 
 a\in \mathcal{V}$ is a MDS containing $a.$ Therefore, $\mathpzc{d}_{\mathpzc{d}}(a)=1.$
\end{proof}
\begin{prop}\label{cmpb}
For any vertex $a\in \mathcal{K}_{n_1,n_2}$ of order $n_1+n_2$, $n_1,n_2>1$, $\mathpzc{d}_{\mathpzc{d}}(a)=2.$
\end{prop}
\begin{proof}
    Let $a\in \mathcal{K}_{n_1,n_2}$. Then either $a\in \mathcal{V}_1$ or $a\in \mathcal{V}_2$, where $\mathcal{V}_1$ and $\mathcal{V}_2$ are partite sets of $\mathcal{K}_{n_1,n_2}.$ Suppose $a\in \mathcal{V}_1$. Then, $\{a,b\}$ where $b$ is any vertex in $\mathcal{V}_2$ is a MDS containing $a$ with least number of vertices. Similarly the case when $a\in \mathcal{V}_2.$ Hence, $\mathpzc{d}_{\mathpzc{d}}(a)=2.$ 
\end{proof}
\begin{prop}
    For a vertex in $\mathcal{K}_{n_1,n_2,...,n_r}$ of order $n_1+n_2+...+n_r$, $n_i>1, i=1,2,...,r,$ $\mathpzc{d}_{\mathpzc{d}}(a)=2.$
\end{prop}
\begin{prop}
    Let $\mathcal{K}_{1,n}$ be a star of order $n+1,$ then 
    $\mathpzc{d}_{\mathpzc{d}}(a')=
     \begin{cases}
     1 \, & \text{ if $a'$ is the center vertex }\\
     n \, & \text{ otherwise. }  
    \end{cases}
    $
\end{prop}
\begin{proof}
    Let $\mathcal{K}_{1,n}$ be a star. Label the center vertex as $a$ and all other vertices as $a_1,a_2,...,a_n.$ Then $\\mathcal{N}[a]=\{a,a_1,a_2,...,a_n\},$ $\mathcal{N}[a_i]=\{a_i,a\},i=1,2,...,n.$ Hence for the center vertex $\{a\}$ is the MDS containing $a$ having the least number of vertices. For any other vertex $a_i,$ since $\mathcal{N}[a_i]\setminus \mathcal{N}[a]=\phi$, $\{a_i,a\}$ is not a MDS. Hence, $\{a_1,a_2,...,a_n\}$ is the MDS containing $a_i$ having the least number of vertices. Therefore, $\mathpzc{d}_{\mathpzc{d}}(a')=
     \begin{cases}
     1 \, & \text{ if $a'$ is the center vertex }\\
     n \, & \text{ otherwise. }  
    \end{cases}
    $
\end{proof}
\begin{cor}
    For a tree of order $n,$ $\lceil{\frac{n}{3}}\rceil\leq \mathpzc{d}_{\mathpzc{d}}(a)\leq (n-1),$ where $a$ is an end vertex. 
\end{cor}
\begin{proof}
    Consider a tree $\mathcal{T}$ on $n$ vertices. Suppose there are $k$ end vertices in $\mathcal{T}$, then at least $k$ vertices are required to dominate all the end vertices. Among all the trees with $n$ vertices a star graph is the one with maximum number of end vertices. The number of end vertices in a star graph with $n$ vertices is $n-1.$ The least number of end vertices in a tree $\mathcal{T}$ is 2, when $\mathcal{T}$ is a path. In a path on $n $ vertices from Theorem \ref{path}, the minimum of  domination degree of an end vertex is $\lceil{\frac{n}{3}}\rceil.$ Hence $\lceil{\frac{n}{3}}\rceil\leq \mathpzc{d}_{\mathpzc{d}}(a)\leq (n-1)$ if $a$ is an end vertex of $\mathcal{T}$. 
\end{proof}
\begin{rem}\label{cycled}
    In [cite chartrand ], it is proved that the domination number of a cycle is $\lceil{\frac{n}{3}}\rceil$ in Example 13.2. 
\end{rem}
\begin{thm}\label{cycle}
Let $\mathcal{C}_n$ be a cycle on $n$ vertices. For any vertex $a \in \mathcal{V}(\mathcal{C}_n)$, $\mathpzc{d}_{\mathpzc{d}}(a)=\lceil{\frac{n}{3}}\rceil.$  
\end{thm}
\begin{proof}
By Remark \ref{cycled} the domination number $\gamma(\mathcal{C}_n)=\lceil{\frac{n}{3}}\rceil.$ Hence for any vertex $a \in \mathcal{V}(\mathcal{C}_n)$, $\mathpzc{d}_{\mathpzc{d}}(a)\geq \lceil{\frac{n}{3}}\rceil.$ Now consider any vertex $a_i,$ $i=1,2,...,n.$ Here, $n=3s+t, $ where $t=0,1,2.$ For the case when $t=0,$ by the proof of \ref{cycled} a minimum dominating set with $\lceil{\frac{n}{3}}\rceil$ is obtained by choosing any vertex and every third vertex beginning at the chosen vertex and moving in a particular direction cyclically. Hence choose $a_i$ as the starting vertex and choose every third vertex moving cyclically in one particular direction. Similarly the case when $t=1,2.$ Choose $a_i$ as the starting vertex and choose every third vertex until $s+1$ vertices are obtained in total. In all the three cases the obtained set is a minimum dominating set containing $a_i$. Hence the obtained set is a MDS. Therefore, $\mathpzc{d}_{\mathpzc{d}}(a_i)=\lceil{\frac{n}{3}}\rceil$ $i=1,2,...,n.$
\end{proof}
\begin{prop}\label{wheel}
    For a wheel graph $\mathcal{W}_n$, 
    $\mathpzc{d}_{\mathpzc{d}}(a')=
     \begin{cases}
     1 \, & \text{ if $a'$ is the center vertex }\\
    \lceil{\frac{n}{3}}\rceil \, & \text{ otherwise. }  
    \end{cases}
    $
\end{prop}
\begin{proof}
Label the center vertex of the wheel graph $\mathcal{W}_n$ as $a$, and the remaining vertices as $a_1,a_2,...,a_n.$ Here, $\mathcal{N}[a]=\{a,a_1,a_2,...,a_n\}$ and every other vertex $a_i$ is adjacent to $a$ and two other vertices as in the case of a cycle. Therefore, $\mathpzc{d}_{\mathpzc{d}}(a)=1.$ Consider vertex $a_i$. Then the set $\{a,a_i\}$ is not a MDS. Hence vertex $a$ cannot be considered while choosing a MDS for vertex $a_i.$ So, the MDS considered in the case of a cycle is taken for vertices of the form $a_i,\, i=1,2,...,n.$ Therefore, $\mathpzc{d}_{\mathpzc{d}}(a_i)=n,\, i=1,2,...,n.$
\end{proof}
\begin{thm}\label{path}
Let $\mathcal{P}_n$ be a path of order $n+1.$ Then, 
$\mathpzc{d}_{\mathpzc{d}}(a_i)=
\begin{cases}
 \lceil{\frac{n}{3}}\rceil \text{ or } \lceil{\frac{n}{3}}\rceil+1 \, &\text{ for } n=3k\\
 \lceil{\frac{n}{3}}\rceil \, &\text{ for } n=3k+1\\
 \lceil{\frac{n}{3}}\rceil \text{ or } \lceil{\frac{n}{3}}\rceil+1 \, &\text{ for } n=3k+2.
 \end{cases}
 $\\
 Moreover, for $n=3k,$  $\mathpzc{d}_{\mathpzc{d}}(a_i)= \lceil{\frac{n}{3}}\rceil$ if $i=3j+2, \, j\leq k.$ And $\mathpzc{d}_{\mathpzc{d}}(a_i)= \lceil{\frac{n}{3}}\rceil+1$, if $i=3j+1 \text{ or } 3j, j\leq k.$\\
 For $n=3k+2,$ $\mathpzc{d}_{\mathpzc{d}}(a_i)= \lceil{\frac{n}{3}}\rceil$ if $i=3j+t, \, t=1,2, j\leq k.$ 
 And $\mathpzc{d}_{\mathpzc{d}}(a_i)= \lceil{\frac{n}{3}}\rceil+1$, if $i=3j, j\leq k.$
 \end{thm}
 \begin{proof}
     Let $\mathcal{P}_n$ be path of order $n.$ From \ref{domin} $\gamma(\mathcal{P}_n)\geq \lceil{\frac{n}{3}}\rceil.$ Hence $\mathpzc{d}_{\mathpzc{d}}(a)\geq \lceil{\frac{n}{3}}\rceil$ for every vertex in the path. Now, label the vertices as $a_1,a_2,...,a_n$ where $a_1$ and $a_n$ are end vertices. Here, the vertices can be categorised into 3: vertices $a_i$ such that $i=3j+t,$ $t=0,1,2.$ 
     For $n$ there are three cases:\\
    \textbf{Case 1}: Order $n=3k.$\\
    First consider the vertices of the form $a_i,\, i=3j+2, \, j\leq k.$ Then $a_2$ dominates $a_1,a_3$ other than $a_2$. Now choose $a_5.$ Then $a_5$ dominates $a_4$ and $a_6$. Continue choosing every third vertex starting from $a_2$. Finally $a_{3k-1}$ is obtained which dominates $a_{3k-2}$ and $a_{3k}$. Hence the obtained set is $\{a_2,a_5,...a_{3k-1}\}$, which is a MDS containing vertices of the form $a_i,\, i=3j+2, \, j\leq k.$ with $k=\lceil{\frac{n}{3}}\rceil$ vertices. Therefore, $\mathpzc{d}_{\mathpzc{d}}(a_i)= \lceil{\frac{n}{3}}\rceil$ if $i=3j+2, \, j\leq k.$ \\
    Now, consider vertices of the form $a_i,\, i=3j+1, \, j\leq k.$ Then $a_1$ dominates $a_2$ other than $a_1$. Now choose $a_4.$ Then $a_4$ dominates $a_3$ and $a_5$. Continue choosing every third vertex starting from $a_1$. Finally $a_{3k-2}$ is obtained which dominates $a_{3k-3}$ and $a_{3k-1}$. The set $\{a_1,a_4,...,a_{3k-2}\}$ contains $k$ vertices, but $a_{3k}$ remains to be dominated. Hence $\{a_1,a_4,...,a_{3k-2},a_{3k}\}$ is a MDS with least number of vertices containing vertices of the form $a_i,\, i=3j+1, \, j\leq k$ with $k+1=\lceil{\frac{n}{3}}\rceil+1$ vertices. Similarly consider vertices of the form $a_i,\, i=3j, \, j\leq k.$ Starting with $a_3.$ Vertex $a_3$ dominates $a_2$ and $a_4$ other than $a_3.$ Continue choosing every third vertex starting from $a_3$. Then vertex $a_6$ dominates $a_5$ and $a_7$. Finally the set $\{a_3,a_6,...,a_{3k}\}$ with $k$ vertices is obtained which dominates every vertex other than $a_1.$ Hence $\{a_1,a_3,a_4,...a_{3k}\}$ is a MDS with the least number of vertices containing vertices of the form $a_i,\, i=3j, \, j\leq k.$ with $k+1=\lceil{\frac{n}{3}}\rceil+1$ vertices. 
    Therefore, $\mathpzc{d}_{\mathpzc{d}}(a_i)= \lceil{\frac{n}{3}}\rceil+1$ if $i=3j+1 \text{ or } 3j, \, j\leq k.$\\
    \textbf{Case 2}: Order $n=3k+1.$\\
    The proof is similar to the proof of case 1. Categorise the vertices into 3 as $a_i,$ $i=3j+s, \, s=0,1,2$ and $j\leq k$, and obtain the MDS in each category. Following the proof as in case 1, the MDS for vertices of the form $a_i,\,i=3j, j\leq k$ with the least number of vertices is $\{a_3,a_6,...a_{3k}\}\cup\{a_1\}$ with $k+1=\lceil{\frac{n}{3}}\rceil$ vertices. The required MDS for vertices of the form $a_i,\,i=3j+1, j\leq k$ is $\{a_1, a_4,..., a_{3k-2}\}\cup\{a_{3k}\}$ with $k+1=\lceil{\frac{n}{3}}\rceil$ vertices. Similarly, the required MDS for vertices of the form $a_i,\,i=3j+2, j\leq k$ is $\{a_2, a_5,..., a_{3k-1}\}\cup\{a_{3k+1}\}$ with $k+1=\lceil{\frac{n}{3}}\rceil$ vertices. Hence for all vertices $\mathpzc{d}_{\mathpzc{d}}(a_i)= \lceil{\frac{n}{3}}\rceil.$ \\
    \textbf{Case 3}: Order $n=3k+2.$ \\
    The proof is similar to the proof of case 1 and 2. The required MDS for vertices of the form $a_i,\,i=3j, j\leq k$ is $\{a_3, a_6,..., a_{3k}\}\cup\{a_1,a_{3k+2}\}$ with $k+2=\lceil{\frac{n}{3}}\rceil+1$ vertices.
    Hence, $\mathpzc{d}_{\mathpzc{d}}(a_i)= \lceil{\frac{n}{3}}\rceil$ if $i=3j, \, j\leq k.$ 
    The required MDS for vertices of the form $a_i,\,i=3j+1, j\leq k$ is $\{a_1, a_4,..., a_{3k-2},a_{3k+1}\}$ with $k+1=\lceil{\frac{n}{3}}\rceil$ vertices. Similarly, the required MDS for vertices of the form $a_i,\,i=3j+2, j\leq k$ is $\{a_2, a_5,..., a_{3k-1}\}\cup \{a_{3k+2}\}$ with $k+1=\lceil{\frac{n}{3}}\rceil$ vertices. Therefore, $\mathpzc{d}_{\mathpzc{d}}(a_i)= \lceil{\frac{n}{3}}\rceil+1$ if $i=3j+1 \text{ or } 3j+2, \, j\leq k.$
 \end{proof}
\begin{thm}\label{book}
Let $\mathcal{B}_n$ be a book graph. Then, 
$\mathpzc{d}_{\mathpzc{d}}(a')=
\begin{cases}
2 \, & \text{ if $a'$ is a center vertex }\\
n \, & \text{ otherwise .} 
\end{cases}
$
\end{thm}
\begin{proof}
Consider the book graph $\mathcal{B}_n= \mathcal{S}_{n+1}\times \mathcal{P}_2.$ Label the vertex of $\mathcal{S}_{n+1}$ with degree $n$ as $a$, and all the remaining vertices as $a_1,a_2,...,a_n.$ Also, label the vertex of $\mathcal{P}_2$ as $b_1,b_2.$ Then in the book graph the center vertex is $ab_1$ and $ab_2.$ Also $ab_1$ is adjacent to $\{a_1b_1,a_2b_1,...,a_nb_1\}.$ Similarly, $ab_2$ is adjacent to $\{a_1b_2,a_2b_2,...,a_nb_2\}.$ So, $ab_1$ and $ab_2$ dominates all vertices of $\mathcal{B}_n.$
Hence, the MDS containing $ab_1$ and $ab_2$, with least number of vertices is $\{ab_1,ab_2\}.$ Therefore, $\mathpzc{d}_{\mathpzc{d}}(ab_1)=\mathpzc{d}_{\mathpzc{d}}(ab_2)=2.$ Now consider a vertex of the form $a_ib_1.$ Then $\mathcal{N}[a_ib_1]=\{a_ib_1,ab_1,a_ib_2\}.$ But $\{a_ib_1,ab_1,ab_2\}$ is not a MDS. Now, MDSs containing $n$ vertices can be obtained by selecting one vertex from each section. Here, $\{a_ib_1,\}\cup\{a_jb_2:\,j\neq i\}$ is a MDS containing $a_ib_1$ with n vertices. All other MDS containing $a_ib_1$ contain $n+1$ vertices. Hence $\{a_ib_1,\}\cup\{a_jb_2:\,j\neq i\}$ is a required MDS with least number of vertices. Hence $\mathpzc{d}_{\mathpzc{d}}(a_ib_1)=n, \,i=1,2,...,n.$ The case is similar for vertices of the form $a_ib_2.$ The set $\{a_ib_2,\}\cup\{a_jb_1:\,j\neq i\}$ is a MDS with least number of vertices containing vertices of the form $a_ib_2.$  Hence $\mathpzc{d}_{\mathpzc{d}}(a_ib_2)=n, \,i=1,2,...,n.$ Therefore $\mathpzc{d}_{\mathpzc{d}}(a')=
\begin{cases}
2 \, & \text{ if $a'$ is a center vertex }\\
n \, & \text{ otherwise .} 
\end{cases}
$ 
\end{proof}
\begin{thm}\label{windmill}
 For a Windmill graph $\mathcal{W}d(r,s)$, 
 $\mathpzc{d}_{\mathpzc{d}}(a')=
\begin{cases}
1 \, & \text{ if $a'$ is the center vertex }\\
s \, & \text{ otherwise .} 
\end{cases}
$
\end{thm}
\begin{proof}
Consider the Windmill graph $\mathcal{W}d(r,s)$. Label the center vertex as $a$ and the vertices in the $i^{th}$ copy as $a_{i_1},a_{i_2},...,a_{i_r}.$ Since the center vertex is adjacent to all other vertices, $\mathpzc{d}_{\mathpzc{d}}(a)=1.$ Consider any other vertex $a_{i_l},i=1,2,...,s \text{ and } l=1,2,...,r.$ Then, $\mathcal{N}[a_{i_l}]=\{a_{i_l}\}\cup\{a_{i_{l'}}: \, l'\neq l\}\cup\{a\}.$ But the set $\{a_{i_l},a\}$ is not a MDS. Hence one vertex from each copies of $\mathcal{K}_r$ is required to dominate all vertices holding the minimality condition. Therefore any MDS with the least number of vertices, containing $a_{i_l}$ contains $s$ vertices. The set $\{a_{i_l}\}\cup\{a_{j_l}: \, j\neq i\}$ is one such set. Hence, $\mathpzc{d}_{\mathpzc{d}}(a_{i_l})=s, i=1,2,...,s \text{ and } l=1,2,...,r.$ Therefore, $\mathpzc{d}_{\mathpzc{d}}(a')=
\begin{cases}
1 \, & \text{ if $a'$ is the center vertex }\\
s \, & \text{ otherwise .} 
\end{cases}
$
\end{proof}
\begin{thm}\label{ktree}
 For a Kragujevac tree $\mathcal{T},$  of order $n=1 +\sum\limits_{i=1}^t (2s_i + 1),$  
 $\mathpzc{d}_{\mathpzc{d}}(a')= s_1+s_2+...+s_t+1
$ where $s_i$ is the number of end vertices in $i^{th}$ branch of $\mathcal{T}, \forall a'\in \mathcal{V}(\mathcal{T}).$ 
\end{thm}
\begin{proof}
    Consider a Kragujevac tree $\mathcal{T},$ of order $n=1 +\sum\limits_{i=1}^t (2s_i + 1).$ Label the central vertex of $\mathcal{T}$ as $a$, then end vertices of the $i^{th}$ branch  as $\{a_{i1},a_{i2},...,a_{is_i}\}.$  The neighbors of the end vertices of the $i^{th}$ branch as $\{b_{i1},b_{i2},...,b_{is_i}\}$ respectively. Also label the $t$ neighbors of $a$ as $c_1,c_2,...,c_t.$ Now consider any end vertex $a_{ij},i=1,2,...,t \text{ and } j\in\{1,2,...,s_i\}.$ Then either $a_{ij}$ or $b_{ij}$ is required to dominate $a_{ij}$. Hence to dominate the end vertices, $k=s_1+s_2+...+s_t$ vertices are required. Consider the branch containing $a_{ij}.$ Choose $\mathcal{A}=a_{ij}\cup \{b_{il'}: \,l'\neq l\}$. Then $\mathcal{A}$ dominates every vertices of branch $i.$ Now for any other branch choose all $b_{i'l}, i'\neq i, l\in\{1,2,...,s_{i'}\}.$ By choosing such a set every vertex except the central vertex is dominated. Hence to dominate every vertex minimum of $k+1$ vertices are required. Hence $\{a_{ij}\}\cup \bigcup\limits_{\substack{j'\in\{1,2,...,s_i\}\\j'\neq j }}b_{ij'}\cup \bigcup\limits_{\substack{l\in\{1,2,...,s_{i'}\}\\i'\neq i }}b_{i'l}\cup \{a\}$ is a MDS containing $a_{ij}$ with least number of vertices. Therefore $\mathpzc{d}_{\mathpzc{d}}(a')= k+1,$ if $a'$ is an end vertex. The case is similar for the vertices of the form $b_{ij}.$ Choose all $b_{ij}, $ $i=1,2,...,t $ and $j\in\{1,2,...,s_i\}$ and $\{a\}$ for a required MDS containing $b_{ij}$ with least number of vertices. Now, for the central vertex any of the above mentioned sets can be considered to get the MDS containing $a.$ Therefore $\mathpzc{d}_{\mathpzc{d}}(a')= s_1+s_2+...+s_t+1
$ where $s_i$ is the number of end vertices in $i^{th}$ branch of $\mathcal{T}, \forall a'\in \mathcal{V}(\mathcal{T}).$   
\end{proof}
\begin{rem}\label{petersen}
 Consider the Petersen graph with labelling as in Figure \ref{petersenf}. The domination number of Petersen graph is 3. Consider vertex $a_1.$ From Proposition \ref{in1}, $\mathpzc{d}_{\mathpzc{d}}(a_1)\geq 3.$ Now $\{a_1,a_7,a_{10}\}$ is a MDS with 3 vertices containing $a_1. $ Hence $\mathpzc{d}_{\mathpzc{d}}(a_1)= 3.$ Similarly, if $a_9$ is considered, $\{a_9,a_1,a_2\}$ is a required MDS. By the same argument $\mathpzc{d}_{\mathpzc{d}}(a_i)= 3, \forall i=1,2,...,10.$ The domination number of Herschel graph and Gr\"{o}tzsch graph is 3. And the domination degree of every vertex in Herschel graph and Gr\"{o}tzsch graph is aslo 3. 
\end{rem}
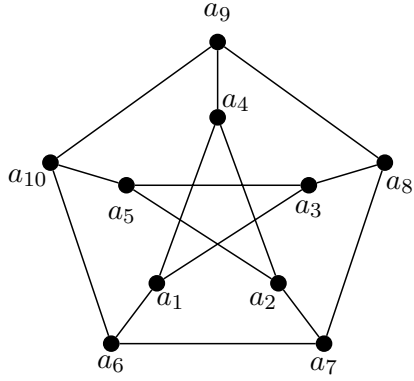
\begin{figure}
\begin{center}
    \begin{tikzpicture}[scale=0.2,inner sep=2pt,line width=0.2mm]
    \draw (0,0) node(1) [circle, draw, fill] {}
      (8,0) node(2) [circle, draw, fill] {}
      (10,6.5) node(3) [circle, draw, fill] {}
      (4,11) node(4) [circle, draw, fill] {}
      (-2,6.5) node(5) [circle, draw, fill] {}
      (-3,-4) node(6) [circle, draw, fill] {}
      (11,-4) node(7) [circle, draw, fill] {}
      (15,8) node(8) [circle, draw, fill] {}
      (4,16) node(9) [circle, draw, fill] {}
      (-7,8) node(10) [circle, draw, fill] {};
    \draw []  (6) to (7) to (8) to (9) to (10) to (6) to (1) to (4) to (2) to (5) to (3) to (1)
              (2) to (7)
              (3) to (8)
              (4) to (9)
              (5) to (10); 
    \draw (0.9,-0.2) node[below] {$a_1$}
          (7,-0.2) node[below] {$a_2$}
          (10,4) node[above] {$a_3$}
          (5.2,11) node[above] {$a_4$}
          (-2.3,5.5) node[below] {$a_5$}
          (-3,-6.3) node[above] {$a_6$}
          (11,-4.5) node[below]{$a_7$}
          (16,7.5) node[below]{$a_8$}
          (4,17) node[above]{$a_9$}
          (-8.5,7) node[]{$a_{10}$}
          ;
\end{tikzpicture}
\caption{Petersen graph}
    \label{petersenf}
\end{center}
\end{figure}
\begin{defn}\label{regular}
 A graph is $k-$ domination regular graph (k-DRG) or simply domination regular graph (DRG), if $_{\mathpzc{d}}\mathpzc{d}(a)=k$, $\forall a\in \mathcal{V}(\mathcal {G})$.    
\end{defn}
\begin{exam}
    An example of DRG is a cycle. By Theorem \ref{cycle} it is clear that for every vertex $a$ $\mathpzc{d}_{\mathpzc{d}}(a)=\lceil{\frac{n}{3}}\rceil.$ Similarly, Petersen graph, Herschel graph and Gr\'{o}tzsch graph are also examples of DRG. 
\end{exam}
\begin{thm}\label{union}
    Let $\mathcal{G}=\bigcup\limits_{i=1}^t \mathcal{G}_i$ be the disjoint union of graphs $\mathcal{G}_1, \mathcal{G}_2,...,\mathcal{G}_t$, then $$_{\mathcal{G}}\mathpzc{d}_{\mathpzc{d}}(a)=_{\mathcal{G}_i}\mathpzc{d}_{\mathpzc{d}}(a)+\sum\limits_{\substack{{j=1}\\ j\neq i}}^t \gamma(\mathcal{G}_j),$$ 
    where $n$ is the order of $\mathcal{G}$ and $n_i$ is order of $\mathcal{G}_i.$
\end{thm}
\begin{proof}
    Let $\mathcal{G}=\bigcup\limits_{i=1}^t \mathcal{G}_i$ be the disjoint union of graphs $\mathcal{G}_1, \mathcal{G}_2,...,\mathcal{G}_t$.
    Consider a vertex $a$ and suppose $a\in \mathcal{V}(\mathcal{G}_k).$ A MDS of $\mathcal{G}$ is disjoint union of MDS of each component $\mathcal{G}_i.$  Let the MDS containing $a$, considered for $_{\mathcal{G}_k}\mathpzc{d}_{\mathpzc{d}}(a)$ be $\mathcal{D}.$ Since the disjoint union of graphs is considered, in any other component $\mathcal{G}_l, l\neq k,$ the MDS contributed for $_{\mathcal{G}}\mathpzc{d}_{\mathpzc{d}}(a)$ will be the minimum dominating set $\mathcal{D}_l$ of $\mathcal{G}_l$. Hence $\mathcal{D}\cup \bigcup\limits_{\substack{{j=i}\\j\neq k}}^t\mathcal{D}_j$ is the required MDS containing $a$ with least number of vertices. Hence $\gamma(\mathcal{G}_l), l\neq k$ gets added to $_{\mathcal{G}_k}\mathpzc{d}_{\mathpzc{d}}(a)$. 
    Therefore $$_{\mathcal{G}}\mathpzc{d}_{\mathpzc{d}}(a)=_{\mathcal{G}_i}\mathpzc{d}_{\mathpzc{d}}(a)+\sum\limits_{\substack{{j=1}\\ j\neq i}}^t \gamma(\mathcal{G}_j).$$
\end{proof}
\begin{thm}
    Consider the join $\mathcal{G}=\mathcal{G}_1+\mathcal{G}_2$ of  $\mathcal{G}_1$ and $\mathcal{G}_2. $ Then
    $_{\mathcal{G}}\mathpzc{d}_{\mathpzc{d}}(a)=
     \begin{cases}
         1 \, &\text{ if } _{\mathcal{G}_i}\mathpzc{d}_{\mathpzc{d}}(a)=1, i=1,2.\\
         2 \, &\text{ otherwise .}
        
     \end{cases}
    $
\end{thm}
\begin{proof}
    Let $\mathcal{G}_1$ and $\mathcal{G}_2$ be two graphs and $\mathcal{G}=\mathcal{G}_1+\mathcal{G}_2.$ Consider a vertex $a$ in $\mathcal{G}.$ Then either $a$ is a vertex of $\mathcal{G}_1$ or $\mathcal{G}_2.$ Suppose $a\in \mathcal{V}(\mathcal{G}_1)$. If $_{\mathcal{G}_1}\mathpzc{d}_{\mathpzc{d}}(a)=1$, it means that $a$ dominates all vertices of $\mathcal{G}_1.$ By the definition of join $a$ is adjacent to all vertices of $\mathcal{G}_2$ also. Hence $_{\mathcal{G}}\mathpzc{d}_{\mathpzc{d}}(a)=1$. The case is similar if $a\in \mathcal{V}(\mathcal{G}_2) .$ Now, $_{\mathcal{G}_1}\mathpzc{d}_{\mathpzc{d}}(a) \neq 1$, then $\{a,b\}, b\in \mathcal{V}(\mathcal{G}_2)$ dominates all vertices of $\mathcal{G}.$ Hence $\{a,b\}$ is a MDS containing $a$ with the least number of vertices. Therefore $_{\mathcal{G}}\mathpzc{d}_{\mathpzc{d}}(a)=2.$ Hence $_{\mathcal{G}}\mathpzc{d}_{\mathpzc{d}}(a)=
     \begin{cases}
         1 \, &\text{ if } _{\mathcal{G}_i}\mathpzc{d}_{\mathpzc{d}}(a)=1, i=1,2.\\
         2 \, &\text{ otherwise .}
        
     \end{cases}
    $
\end{proof}
\begin{thm}
    Let $\mathcal{G}$ be a graph. Consider the composition $\mathcal{H}=\mathcal{G}\circ \mathcal{K}_n.$ Let the vertices of $\mathcal{G}$ be $a_1,a_2,...,a_m$ and the vertices of $\mathcal{K}_n$ be $b_1,b_2,...,b_n.$  Then, 
    $$_{\mathcal{H}}\mathpzc{d}_{\mathpzc{d}}(a_i,b_j)= _{\mathcal{G}}\mathpzc{d}_{\mathpzc{d}}(a_i).$$
\end{thm}
\begin{proof}
    Consider the composition graph $\mathcal{H}.$ Any vertex of $\mathcal{H}$ is of the form $(a_i,b_j), i=1,2,...,m$ and $j=1,2,...,n.$ By the definition of composition a vertex $(a_i,b_j) $ dominated all vertices of the form $(a_i,b_{j'}), j'=1,2,...,n$ and $(a_{i'},b_{j'}), j'=1,2,...,n,$ $a_{i'}$ are such that $a_i$ is adjacent to $a_{i'}$ in $\mathcal{G}.$ Hence any vertex $(a_i,b_j)$ is not adjacent to a vertex $(a,b)$ if $a_i$ is not adjacent to $a.$ Suppose that a MDS containing $a_i$ considered for $_{\mathcal{G}}\mathpzc{d}_{\mathpzc{d}}(a_i)$ is $\{a_i,c_1,c_2,...c_k\}, c_1,c_2,...,c_k\in \{a_1,a_2,...,a_m\}\setminus \{a_i\}.$ Then $\mathcal{D}=$ $\{(a_i,b_j),(c_1,b_1),...,(c_k,b_1)\}$ is a DS containing $(a_i,bj).$ Also $\mathcal{D}$ is minimal. Since if $\mathcal{D}\setminus \{(c_{k'},b_1)\}$ is a DS, then $\{a_i,c_1,...c_k\}\setminus \{c_{k'}\}$ is also a DS which is not possible. By the similar argument if $\mathcal{D}'$ is a MDS containing $(a_i,b_j)$ having less number vertices than $\mathcal{D}$, then collecting all the first coordinate vertices of $\mathcal{D}'$ a MDS of $\mathcal{G}$ is obtained containing $a_i$ with less number of vertices than $\{a_i,c_1,c_2,...c_k\}$ which is also not possible. Hence $\mathcal{D}=$ $\{(a_i,b_j),(c_1,b_1),...,(c_k,b_1)\}$ is a required set. Therefore $$_{\mathcal{H}}\mathpzc{d}_{\mathpzc{d}}(a_i,b_j)= _{\mathcal{G}}\mathpzc{d}_{\mathpzc{d}}(a_i).$$
\end{proof}
\begin{thm}
    Let $\mathcal{H}=\mathcal{G}_1\odot \mathcal{G}_2$ be the corona of two graphs $\mathcal{G}_1$ and $\mathcal{G}_2.$ Let the vertices of $\mathcal{G}_1$ be $a_1,a_2,...,a_n$ and the vertices of $\mathcal{G}_2$ be $b_1,b_2,...,b_m $. In $\mathcal{H}$, there are $n$ copies of $\mathcal{G}_2$. Let the vertices in the $i^{th}$ copy be $b_{i1},b_{i2},...,b_{im}$. Then, 
    $$_{\mathcal{G}}\mathpzc{d}_{\mathpzc{d}}(a)=
     \begin{cases}
         n \, &\text{ if } a=a_i, i=1,2,...,n\\
         _{\mathcal{G}_2}\mathpzc{d}_{\mathpzc{d}}(a)+(n-1) \, &\text{ if } a=b_{ij}, i=1,2,...,n \text{ and } j=1,2,...,m.
     \end{cases}
    $$
\end{thm}
\begin{proof}
    Let $\mathcal{H}=\mathcal{G}_1\odot \mathcal{G}_2$ be the corona of two graphs $\mathcal{G}_1$ and $\mathcal{G}_2.$ Any vertex of the form $a_i,i=1,2,...,n$ dominates all vertices in the $i^{th}$ copy and also the vertices which are adjacent to $a_i$ in $\mathcal{G}_1.$ To dominate any vertex of the form $b_{ij}$ either vertex $a_i$ or any other vertex of the form $b_{ij'}$ such that $b_j$ is adjacent to $b_{j'}$ in $\mathcal{G}_2$ is required. This means that to dominate all vertices of the form $b_{ij}$ at least $n$ vertices are required. Hence, if $a_i$ is considered, a MDS containing $a_i$ with the least number of vertices is $\{a_1,a_2,...,a_n\}$. Now, consider any $b_{ij}.$ To dominate all other vertices of the $i^{th}$ copy either $a_i$ or any vertex $b_{ij'}$ such that $b_j$ dominates $b_{j'}$ in $\mathcal{G}_2$ is needed. Also by choosing all other $a_{i'}, i'\neq i$, the remaining $b_{i'j'}, j'=1,2,...,m$ are dominated with $n-1$ vertices which is the minimum. But $\{b_{ij}\}\cup \bigcup\limits_{i\neq i'}a_{i'}$ is not a MDS. Hence the dominating set considered for $_{\mathcal{G}_2}\mathpzc{d}_{\mathpzc{d}}(b_j)$ in $\mathcal{G}_2$ is required to dominate the remaining $b_{ij'}, j'\neq j$  in the $i^{th}$ copy. Therefore $_{\mathcal{H}}\mathpzc{d}_{\mathpzc{d}}(b_{ij})=_{\mathcal{G}_2}\mathpzc{d}_{\mathpzc{d}}(b_j)+(n-1).$ Hence, $_{\mathcal{G}}\mathpzc{d}_{\mathpzc{d}}(a)=
     \begin{cases}
         n \, &\text{ if } a=a_i, i=1,2,...,n\\
         _{\mathcal{G}_2}\mathpzc{d}_{\mathpzc{d}}(a)+(n-1) \, &\text{ if } a=b_{ij}, i=1,2,...,n \text{ and } j=1,2,...,m.
     \end{cases}
    $
\end{proof}
\section{Domination index in graphs}
The section discusses new concept of domination index of a graph using the domination degree of a vertex. Various bounds for the $DI$ using existing graph parameters are discussed. Using the results in Section 3 the $DI$ of various significant graphs are studied. The study is also conducted on spanning subgraphs, spanning trees of a graph. The operations on graphs are also considered in the study.  
\begin{defn}\label{di}
Let $\mathcal{G}$ be a graph. Then domination index (DI) of $\mathcal{G}$, denoted as $DI(\mathcal{G})$ is defined as the sum of domination degree of vertices of $\mathcal{G},$ i.e,
$$DI(\mathcal{G})=\sum\limits_{a\in \mathcal{V}(\mathcal{G})} \mathpzc{d}_{\mathpzc{d}}(a).$$    
\end{defn}
An illustration of Definition \ref{di} is provided in Example \ref{die}.
\begin{exam}\label{die}
For the graph in Figure \ref{f1}, $\mathpzc{d}_{\mathpzc{d}}(a_1)=\mathpzc{d}_{\mathpzc{d}}(a_2)=\mathpzc{d}_{\mathpzc{d}}(a_9)=\mathpzc{d}_{\mathpzc{d}}(a_{10})=\mathpzc{d}_{\mathpzc{d}}(a_7)=\mathpzc{d}_{\mathpzc{d}}(a_8)=3$  and $\mathpzc{d}_{\mathpzc{d}}(a_3)=\mathpzc{d}_{\mathpzc{d}}(a_4)=\mathpzc{d}_{\mathpzc{d}}(a_5)=\mathpzc{d}_{\mathpzc{d}}(a_{6})=4.$ Therefore, $DI(\mathcal{G})=\sum\limits_{a_i\in \mathcal{V}(\mathcal{G})} \mathpzc{d}_{\mathpzc{d}}(a_i)=6\times 3+4\times 4=34.$
\end{exam}
From \ref{in1}, \ref{in2}, \ref{irin} the inequalities for $DI$ follows:
\begin{prop}
    Let $\mathcal{G}$ be a graph of order $n.$ Then
    $\gamma(\mathcal{G})\leq \frac{DI(\mathcal{G})}{n}\leq \Gamma(\mathcal{G}).$
\end{prop}
\begin{prop}
    Let $\mathcal{G}$ be a graph of order $n.$ Then
    $DI(\mathcal{G})\leq n \times WI(\mathcal{G}).$
\end{prop}
\begin{prop}
    Let $\mathcal{G}$ be a graph of order $n.$ Then
    $ir\leq \gamma \leq \frac{DI(\mathcal{G})}{n}\leq \Gamma \leq IR$
\end{prop}
\begin{thm}\label{isomorphism}
Let $\mathcal{G}_1\cong \mathcal{G}_2,$ then $DI(\mathcal{G}_1)=DI(\mathcal{G}_2).$    
\end{thm}
\begin{proof}
    Let $\mathcal{G}_1=(\mathcal{V}_1,\mathcal{E}_1)$ and $\mathcal{G}_2=(\mathcal{V}_2,\mathcal{E}_2)$ be isomorphic and $\theta$ be the bijection from ${\mathcal{V}_1}$ to ${\mathcal{V}_2}$ such that $ab\in \mathcal{E}_1$ iff $\theta(a)\theta(b)\in \mathcal{E}_2$. Since $\mathcal{G}_1$ and $\mathcal{G}_2$ are isomorphic the number of edges incident at $a$ and $\Theta(a)$ are same. Hence if $\{a,a_1,a_2,...,a_n\}$ is a DS containing in $\mathcal{G}_1, $ then $\{\theta(a),\theta(a_1),\theta(a_2),...,\theta(a_n)\}$ is a DS containing $\theta(a)$ in $\mathcal{G}_2.$  Therefore, $_{sd}d(u)=_{sd}d(\theta(u))$. Now,
    \begin{align*}
     DI(\mathcal{G}_1)&=\sum\limits_{a\in V(\mathcal{G}_1)} {\mathpzc{d}_{\mathpzc{d}}(a)} \\
     &= \sum\limits_{\theta(a)\in V(\mathcal{G}2)} {\mathpzc{d}_{\mathpzc{d}}(\theta(a))}\\
     &=DI(\mathcal{G}_2)
    \end{align*}
     Therefore, $DI(\mathcal{G}_1)=DI(\mathcal{G}_2).$
\end{proof}
\begin{prop}
    For $\mathcal{K}_n$ of order $n,$ $DI(\mathcal{K}_n)=n.$
\end{prop}
\begin{prop}
    For $\mathcal{K}_{n_1,n_2,...,n_r}$ of order $n_1+n_2+...+n_r,$ $DI(\mathcal{K}_{n_1,n_2,...,n_r})=2(n_1+n_2+...+n_r).$
\end{prop}
\begin{prop}
    For a star graph $\mathcal{K}_{1,n},$ $DI(\mathcal{K}_{1,n})=1+n^2.$
\end{prop}
\begin{prop}
    For a cycle $\mathcal{C}_n,$ $DI(\mathcal{C}_{n})=n \lceil{\frac{n}{3}}\rceil.$
\end{prop}
\begin{prop}
    For a wheel graph $\mathcal{W}_n,$ $DI(\mathcal{W}_{n})=1+n \lceil{\frac{n}{3}}\rceil.$
\end{prop}
\begin{prop}
    Let $\mathcal{P}_n$ be a path of order $n+1.$ Then, 
$$DI(\mathcal{P}_n)=
\begin{cases}
 k(3k+2) \, &\text{ for } n=3k\\
 (3k+1)(k+1) \, &\text{ for } n=3k+1\\
 k(k+1)+2(k+1)^2 \, &\text{ for } n=3k+2.
 \end{cases}
 $$
\end{prop}
\begin{proof}
    Consider path $\mathcal{P}_n$ and label the vertices as in proof of Theorem \ref{path}. From Theorem \ref{path}, for $n=3k,$  $\mathpzc{d}_{\mathpzc{d}}(a_i)= \lceil{\frac{n}{3}}\rceil$ if $i=3j+2, \, j\leq k.$ And $\mathpzc{d}_{\mathpzc{d}}(a_i)= \lceil{\frac{n}{3}}\rceil+1$, if $i=3j+1 \text{ or } 3j, j\leq k.$\\
   For $n=3k+2,$ $\mathpzc{d}_{\mathpzc{d}}(a_i)= \lceil{\frac{n}{3}}\rceil$ if $i=3j+t, \, t=1,2, j\leq k.$ 
   And $\mathpzc{d}_{\mathpzc{d}}(a_i)= \lceil{\frac{n}{3}}\rceil+1$, if $i=3j, j\leq k.$ \\
   And for $n=3k+1,$ $\mathpzc{d}_{\mathpzc{d}}(a_i)= \lceil{\frac{n}{3}}\rceil$ $\forall i=1,2,...,n.$\\ Now, for $n=3k,$ exactly $k$ vertices have domination degree $\lceil{\frac{n}{3}}\rceil=k,$ and $3k-k=2k$ vertices have domination degree $\lceil{\frac{n}{3}}\rceil+1=k+1.$ Therefore $DI(\mathcal{P}_n)= k\times k + 2k\times (k+1)=3k^2+2k=k(3k+2).$ \\
   For $n=3k+1,$ every vertex have domination degree $\lceil{\frac{n}{3}}\rceil=k+1.$ Therefore $DI(\mathcal{P}_n)=(3k+1)(k+1).$\\
   Similarly, for $n=3k+2,$ exactly $k$ vertices have domination degree $\lceil{\frac{n}{3}}\rceil+1=k+2,$ and $(3k+2)-k=(2k+2)$ vertices have domination degree $\lceil{\frac{n}{3}}\rceil=k+1.$ Therefore $DI(\mathcal{P}_n)=k(k+2)+2(k+1)(k+1)=k(k+2)+2(k+1)^2.$ Hence $$DI(\mathcal{P}_n)=
\begin{cases}
 k(3k+2) \, &\text{ for } n=3k\\
 (3k+1)(k+1) \, &\text{ for } n=3k+1\\
 k(k+1)+2(k+1)^2 \, &\text{ for } n=3k+2.
 \end{cases}
 $$ 
\end{proof}
\begin{prop}
 For a book graph $\mathcal{B}_n,$ $DI(\mathcal{P}_n)=2(n^2+2).$  
\end{prop}
\begin{proof}
    For a book graph $\mathcal{B}_n$ labelled as in Theorem \ref{book}, exactly two center vertices have domination degree 2 and $2(n+1)-2=2n$ vertices have domination degree $n.$ Therefore, $DI(\mathcal{B}_n)=2\times 2+2n\times n=4+2n^2=2(2n+2).$ 
\end{proof}
\begin{prop}
 For a Windmill graph $\mathcal{W}d(r,s),$ $DI(\mathcal{W}d(r,s)=1+rs^2.$   
\end{prop}
\begin{proof}
    The proof follows from Theorem \ref{windmill} and the fact that there are $s$ copies of $\mathcal{K}_r.$ Hence there is one center vertex and $rs$ other vertices.
\end{proof}
\begin{prop}
    For a Kragujevac tree $\mathcal{T},$ $DI(\mathcal{T})=(1 +\sum\limits_{i=1}^t (2s_i + 1))(1+s_1+s_2+...+s_t).$
\end{prop}
\begin{proof}
    Consider a Kragujevac tree $\mathcal{T}$ with $t$ branches. The number of end vertices in the $i^{th}$ branch is $s_i.$ Then by Theorem \ref{ktree}, every vertex have domination degree $1+s_1+s_2+...+s_t.$ The total number of vertices in $\mathcal{T}$ is $1 +\sum\limits_{i=1}^t (2s_i + 1).$ Hence, $DI(\mathcal{T})=(1 +\sum\limits_{i=1}^t (2s_i + 1))(1+s_1+s_2+...+s_t).$
\end{proof}
\begin{rem}
    The $DI$ of Petersen graph is $3\times 10=30.$ Similarly, both Herschel and Gr\"{o}tzsh graph have $DI=33.$ 
\end{rem}
\begin{prop}
    Let $\mathcal{H}$ be a spanning subgraph of $\mathcal{G}$, then $DI(\mathcal{G})\leq DI(\mathcal{H}).$
\end{prop}
\begin{cor}
   Let $\mathcal{T}$ be a spanning tree of $\mathcal{G}$, then $DI(\mathcal{G})\leq DI(\mathcal{T}).$ 
\end{cor}
Let $\mathcal{G}$ and $\mathcal{G}$ be two graphs. Since $\mathcal{G}\times \mathcal{H}\leq \mathcal{G}\boxtimes \mathcal{H}\leq \mathcal{G}\circ\mathcal{H}, $ and $\mathcal{G}\square \mathcal{H}\leq \mathcal{G}\boxtimes \mathcal{H}\leq \mathcal{G}\circ\mathcal{H},$ the result follows:
\begin{prop}
    Let $\mathcal{G}$ and $\mathcal{G}$ be two graphs $ DI(\mathcal{G}\circ\mathcal{H})\leq DI(\mathcal{G}\boxtimes \mathcal{H}) \leq \min\{DI(\mathcal{G}\times \mathcal{H}),DI(\mathcal{G}\square \mathcal{H})\}.$
\end{prop}
\section{Application}
The selection of suitable locations for facilities to meet demand requirements has garnered significant attention among researchers. Within the banking industry, one particular service facility holds immense importance: automatic teller machines (ATMs). These machines offer a wide array of banking services, including balance inquiries, withdrawals, statement requests, and deposits. By providing customers with access to banking services at any time, without the necessity of visiting a nearby bank branch, ATMs greatly enhance convenience. Furthermore, strategically placing ATMs in convenient locations can alleviate the workload of bank branches. Therefore, banks should allocate ATMs in favorable locations, taking into consideration customer concerns. The problem is to select locations to set up ATMs and its bank branch in a city that cater to customer convenience and demand. While selecting location for the bank factors like population, commercial centers, residential regions should be considered.  \\
\begin{figure}[h!]
    \centering
    \includegraphics[height=7cm,width=10cm]{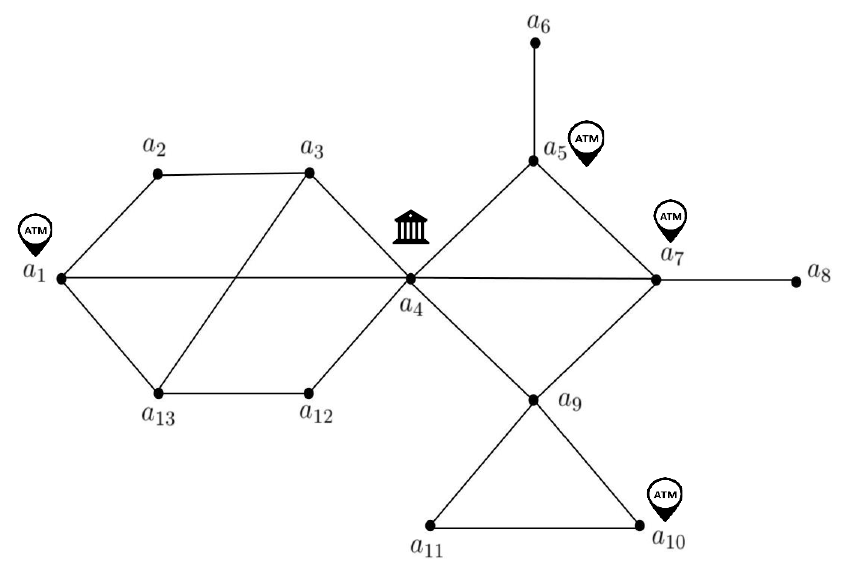}
    \caption{Bank and ATM allocation in a city.}
    \label{app}
\end{figure}
Consider a city and the road network connecting different parts of the city. A graphical representation of the city and road network is used for precision. Vertices of the graph represent different parts of the city and edges represent the roads connecting different places. In a city there often exists a prime location that encompasses a blend of residential areas, commercial hubs, and a constant influx of people. Choose one location that stands out with the highest population and crowd. The idea is to set up the main bank branch in such a location. The ATMs can be placed in such a way that the other parts of the city and the parts where ATMs are located are at most one block away. The notion of domination degree is used to set up the bank and its ATMs in different parts of the city to ensure that all parts are covered effectively. Consider the graph in Figure \ref{app}. Suppose that the part corresponding vertex $a_4$ is the most prominent part of the city. Then the main branch of the bank can be set up at $a_4$, which is most convenient place for the public. Now consider a MDS with least cardinality containing $a_4.$ Then the ATMs can be set up at the parts corresponding to the vertices in the MDS other than $a_4$, so that the entire city is covered. In Figure \ref{app} a MDS containing $a_4$ with least number of vertices is $\mathcal{D}=\{a_1,a_4,a_{10},a_5,a_7\}.$ Hence the main branch can be set up at $a_4$ and the ATMs can be set up at parts corresponding to $a_1,a_{10},a_5$ and $a_7$ which cover the entire city and demand points.
\vspace{0.5cm}

By the same argument the concept of domination degree can be applied for allocating any facility in any context. It can also be used in security allocation in computer networks. 
Identify critical nodes, such as network gateways, authentication servers, or repositories of sensitive data, which, if compromised, would significantly impact the network's security. Identify the node with most confidential and sensitive data and consider a MDS containing vertex corresponding to the most critical node. Then provide high protection to the most critical node and some security measures and resources in the nodes corresponding to other vertices in the MDS. 
 \vspace{1cm}
 
 \textbf{Algorithm to find MDS containing a particular vertex.}
 \\
An algorithm to find a MDS containing a particular vertex is provided and illustrated here. Consider a graph $\mathcal{G}$ of order $n,$ and label the vertices as $a_1,a_2,...,a_n.$ Take any vertex say $a_k.$ Suppose it is required to find a MDS containing $a_k$. Then the following algorithm can be used find the MDS.
\\
Follow the following steps to obtain the required MDS.
\begin{enumerate}
    \item Take vertex $a_{k_0}$, $a_{k_0}\in \{a_1,a_2,...,a_n\}$. If $N[a_{k_0}]=\mathcal{V}$, then $\{a_{k_0}\}$ is a required set.
    \item Suppose $N[a_{k_0}]\neq \mathcal{V}$, take any vertex $a_{k_1}\neq a_{k_0}$ such that \\
    $$N[a_{k_0}]\setminus N[a_{k_1}] \neq \phi \text{ and } N[a_{k_1}]\setminus N[a_{k_0}] \neq \phi$$
    \item If $N[a_{k_1}]\cup N[a_{k_0}]=\mathcal{V}$, then $\{a_{k_1}, a_{k_0}\}$ is a required set. 
    \item Suppose $N[a_{k_1}]\cup N[a_{k_0}]\neq \mathcal{V}$, take next vertex $a_{k_2}\neq a_{k_0}, a_{k_1}$, such that \\
    $$N[a_{k_0}]\setminus\{ N[a_{k_1}] \cup N[a_{k_2}]\}\neq \phi \text{ and } N[a_{k_1}]\setminus \{N[a_{k_2}]\cup N[a_{k_0}]\} \neq \phi \text{ and } N[a_{k_2}]\setminus\{ N[a_{k_1}] \cup N[a_{k_0}]\}\neq \phi .$$
    \item If $N[a_{k_2}] \cup N[a_{k_1}]\cup N[a_{k_0}]=\mathcal{V}$, then $\{a_{k_2},a_{k_1},a_{k_0}\}$ is the required set.
    \item Continue the process until the whole vertex set is obtained. 
\end{enumerate}
Check all the possible combinations of vertices satisfying the conditions in each step of the algorithm. \\
After obtaining the MDS, the cardinality of the set can also be found. \\
After obtaining all the MDS containing $a_{k_0}$ check the cardinality of each set and take the set with least cardinality. The minimum cardinality of MDSs containing $a_{k_0}$ is the domination degree of $a_{k_0}.$

\vspace{0.7cm}
\textbf{Illustration of Algorithm}
\vspace{0.7cm}

\begin{figure}
    \centering
    \includegraphics[height=4cm,width=11cm]{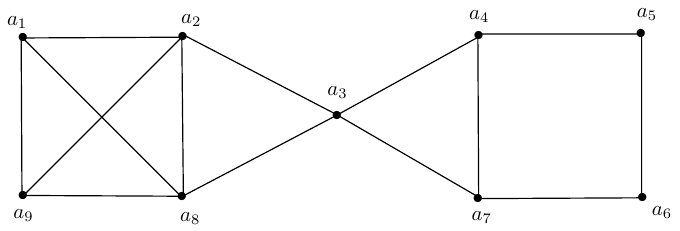}
    \caption{Illustration of the algorithm to find MDS}
    \label{alg}
\end{figure}

Consider the graph in Figure \ref{alg}. In order to find a MDS containing $a_2$ follow the following steps:
First the closed neighborhood of each vertex is obtained. 
From the matrix, $N[a_1]=\{a_1,a_2,a_9,a_8\}$, $N[a_2]=\{a_2,a_1,a_9,a_8,a_3\}$, $N[a_3]=\{a_3,a_2,a_8,a_4,a_7\}$, $N[a_4]=\{a_4,a_3,a_7,a_5\}$, $N[a_5]=\{a_5,a_4,a_6\}$, $N[a_6]=\{a_6,a_5,a_7\}$, $N[a_7]=\{a_7,a_3,a_4,a_6\}$, $N[a_8]=\{a_8,a_3,a_1,a_2,a_9\}$ and $N[a_9]=\{a_9,a_1,a_2,a_8\}.$ 
\begin{enumerate}
    \item Here $N[a_2]\neq \mathcal{V}. $
    \item Take next vertex $a_1.$ But, $N[a_1]\setminus N[a_2]=\phi$.Therefore $a_1$ cannot be considered.
    \item By the same argument $a_9$ and $a_8$ is not considered.
    \item Now, take $a_3$. Here, $N[a_2]\setminus N[a_3]\neq \phi $ and $N[a_3]\setminus N[a_2]\neq \phi $. But $N[a_2]\cup N[a_3]\neq \mathcal{V}$.
    \item Hence consider next vertex $a_4$. Then, $N[a_3]\setminus \{N[a_4]\cup N[a_2]\}= 
    \phi .$ Therefore $a_4$ cannot be considered. By the same argument $a_7$ is not considered. 
    \item Now, consider $a_5,$ then $N[a_5]\setminus \{N[a_2] \cup N[a_3]\}\neq \phi ,$ $N[a_3]\setminus \{N[a_2] \cup N[a_5]\}\neq \phi $ and $N[a_2]\setminus \{N[a_3] \cup N[a_5]\}\neq \phi. $
     Also, $N[a_2]\cup N[a_3]\cup N[a_5]=\mathcal{V}.$
    \item Hence $\{a_2,a_3,a_5\}$ is a required set containing $a_2$ with cardinality 3.
\end{enumerate}
Similarly, $\{a_2, a_4, a_6\}$ is also a set. By, considering all MDSs containing $a_4$ the domination degree of $a_4$ is obtained as 3. Hence the given algorithm can be used to find a minimal MDS containing a particular vertex of a graph. 
\section{Conclusion}
A novel notion of domination degree is introduces in graphs. The idea is used to obtain domination index of a graph. The cardinality of MDS containing a particular vertex is used to define the concepts. The study is conducted on several graph classes like complete graphs, complete bipartite, $r-$ partite graphs, wheels, cycles, paths, book graph, windmill graph, Kragujevac trees. The union, join, composition and corona of graphs are also considered in the study. An application of domination degree of a vertex is discussed in facility allocation problem. An algorithm to find a MDS containing a particular vertex is also provided in the article. 

\vspace{2cm}
\textbf{\large{Acknowledgement}}\\\\
The first author gratefully acknowledges the financial support of Council of Science and Industrial Research (CSIR), Government of India.\\
The authors would like to thank the DST, Government of India, for providing support to carry out this work under the scheme 'FIST' (No.SR/FST/MS-I/2019/40).

\end{document}